  \documentclass[11pt,a4paper,twoside]{article}
  \usepackage[a4paper, hmargin=3.1cm,vmargin=3.82cm]{geometry}

  \usepackage{amsmath,amssymb,amsthm,amsfonts}
  \usepackage{mathabx}
    \usepackage{xspace,xcolor}
  \usepackage{parskip,fancyhdr,url}
  
  \usepackage{fourier}
  \usepackage{enumitem}
    \usepackage{bm}
  \usepackage{todonotes}
  \usepackage[
     breaklinks,colorlinks,
     citecolor=blue,linkcolor=blue,urlcolor=teal,
     pdftitle={Categorical characterisations of quasi-isometric embeddings},pdfauthor={Robert Tang}
  ]{hyperref}
  \usepackage{tikz, tikz-cd}
  \usetikzlibrary{
    arrows, calc, matrix, decorations.pathmorphing, shapes, arrows
    }
  \usepackage{mathtools}
  \usepackage{centernot}

  \setenumerate[0]{label=(\arabic*)}

  \makeatletter
    \def\tagform@#1{\maketag@@@{%
     \textbf{(\ignorespaces#1\unskip\@@italiccorr)}}}%
     \renewcommand{\eqref}[1]{\textup{\maketag@@@{(\ignorespaces%
          {\ref{#1}}\unskip\@@italiccorr)}}}
  \makeatother

  \newcommand\address[1]{}
  \newcommand\email[1]{}
  \newcommand\dedicatory[1]{}

  \theoremstyle{plain}
  \newtheorem{theorem}[equation]{Theorem}
  \newtheorem{proposition}[equation]{Proposition}
  \newtheorem{corollary}[equation]{Corollary}
  \newtheorem{lemma}[equation]{Lemma}

  \theoremstyle{definition}
  \newtheorem{definition}[equation]{Definition}
  \newtheorem{remark}[equation]{Remark}
  \newtheorem{example}[equation]{Example}

  \numberwithin{equation}{section}

  \newcommand{\bsf}[1]{\ensuremath{\bm{\mathsf{#1}}}\xspace}

  \newcommand{\cat}{\ensuremath{\bsf{C}\xspace}}

  \newcommand{\CLip}{\ensuremath{\bsf{CLip}\xspace}}

  \newcommand{\Set}{\ensuremath{\bsf{Set}\xspace}}

  \newcommand{\Born}{\ensuremath{\bsf{Coarse}\xspace}}
  \newcommand{\Top}{\ensuremath{\bsf{Top}\xspace}}

  \newcommand{\barcat}{\overline{\cat}}
  \newcommand{\barCL}{\overline{\CLip}}
  \newcommand{\barBorn}{\overline{\Born}}

   \DeclareMathOperator{\Mor}{Mor} 

  \newcommand{\colimit}{\ensuremath{\varinjlim\xspace}}

      \DeclareMathOperator{\Eq}{Eq} 

  \newcommand{\from}{\ensuremath{\colon \thinspace}\xspace}

  \newcommand{\close}{\ensuremath{\approx}\xspace}

   \DeclareMathOperator{\diam}{diam}

  \newcommand{\N}{\ensuremath{\mathbb{N}}\xspace}
  \newcommand{\R}{\ensuremath{\mathbb{R}}\xspace}

  \newcommand{\Hy}{\ensuremath{\mathbb{H}}\xspace}

  \newcommand{\param}{{\mathchoice{\mkern1mu\mbox{\raise2.2pt\hbox{$
  \centerdot$}}
  \mkern1mu}{\mkern1mu\mbox{\raise2.2pt\hbox{$\centerdot$}}\mkern1mu}{
  \mkern1.5mu\centerdot\mkern1.5mu}{\mkern1.5mu\centerdot\mkern1.5mu}}}

  \fancypagestyle{plain}

  \AtEndDocument{\bigskip{\footnotesize%
  Department of Pure Mathematics, Xi'an Jiaotong--Liverpool University, 111 Ren'ai Road, Suzhou Industrial Park, Suzhou, Jiangsu, 215123, P. R.~China
  \par \texttt{robert.tang@xjtlu.edu.cn}}}

\begin{document}

      \title{Categorical characterisations of quasi-isometric embeddings}
    \author{Robert Tang} \date{}

  \maketitle \thispagestyle{empty}

  \begin{abstract}

  We characterise the (closeness classes of) quasi-isometric embeddings as the regular monomorphisms in the coarsely Lipschitz category, formalising the notion that they are isomorphisms onto their image. Furthermore, we prove that the coarsely Lipschitz category is coregular, and hence admits an (Epi, RegMono)--orthogonal factorisation system. Consequently, quasi-isometric embeddings are equivalently characterised as the effective, strong, or extremal monomorphisms. Finally, we prove that the coarsely Lipschitz category is not coexact in the sense of Barr.

  \end{abstract}

  \section{Introduction}

  The notion of quasi-isometry is ubiquitious throughout geometric group theory and coarse geometry, providing a natural equivalence relation for studying the large-scale geometry of finitely generated groups via their word metrics. Following the foundational work of Gromov in the 1980's \cite{Gro87}, much attention has been devoted to studying quasi-isometry invariants of groups such as hyperbolicity, boundary structure, growth, and amenability \cite{Gro93, BH99, Har00, DK18}. In particular, understanding how these invariants behave under quasi-isometric embeddings is of fundamental importance.

  In this paper, we consider quasi-isometric embeddings in the context of the coarsely Lipschitz category; here, objects are (extended) metric spaces and morphisms are closeness classes of maps admitting an affine upper control.
  A standard result is that (closeness classes of) quasi-isometries are precisely the isomorphisms in this category. Quasi-isometric embeddings are necessarily monomorphisms, however, the converse does not hold. We give some examples in Section \ref{sec:control} which demonstrate this. Our first main result is a characterisation of the quasi-isometric embeddings as the \emph{regular monomorphisms}. In any category, a regular monomorphism is the equaliser of some parallel pair.

  \begin{theorem}[QI-embedding = RegMono]\label{thm:RegMono}
   Let $f \from X \to Y$ be a coarsely Lipschitz map between metric spaces. Then $f$ is a quasi-isometric embedding if and only if its closeness class $\bar f$ is a regular monomorphism in the coarsely Lipschitz category $\barCL$.
  \end{theorem}

    Regular monomorphisms capture the notion of a morphism being an isomorphism onto its image. For example, in the category $\Top$ of topological spaces, the regular monomorphisms are the topological embeddings, whereas the monomorphisms are the continuous injections. Since $\barCL$ is not a concrete category, the appropriate notion of image is not immediately obvious. A key step is to show that every closeness class $\bar f$ in $\barCL$ admits a \emph{regular image}, defined in terms of a universal property (Proposition \ref{prop:image}). A morphism is then a regular monomorphism if and only it factorises through its regular image via an isomorphism (Theorem \ref{thm:regmono}).

  We show further that the factorisation of any morphism through its regular image is canonical. To do so, we show that $\barCL$ is coregular.
  A \emph{coregular category} is a finitely cocomplete category in which all regular images exist and where regular monomorphisms are stable under pushouts. Coregular categories are a generalisation of the abelian categories, with the regular image factorisations playing an analogous role to short coexact sequences in the non-additive setting \cite{BGO71, Bor94, Gra21}.
  Examples of coregular categories include $\Set$ and $\Top$.

  \begin{theorem}[Coarsely Lipschitz category is coregular]\label{thm:Coregular}
  The category $\barCL$ is coregular.
  \end{theorem}

  In any coregular category, the subclasses of epimorphisms and regular monomorphisms together yield an \emph{orthogonal factorisation system} (Definition \ref{def:OFS}). A morphism in $\barCL$ is an epimorphism if and only if any representative is coarsely surjective.

  \begin{corollary}[Factorisation]\label{cor:Orthogonal}
   The category $\barCL$ admits an (Epi, RegMono)--orthogonal factorisation system. In particular, each morphism $\bar f$ admits a factorisation ${\bar f = \bar m \bar e}$, unique up to unique isomorphism, via an epimorphism $\bar e$ and regular monomorphism $\bar m$. \qed
  \end{corollary}

  Another consequence of coregularity is that several subclasses of monomorphisms coincide with the regular monomorphisms. This yields further equivalent characterisations of quasi-isometric embeddings in $\barCL$; see Proposition \ref{prop:equiv_mono} for definitions.

  \begin{corollary}[Equivalent subclasses of monomorphisms]\label{cor:Mono}
   In the category $\barCL$, the following subclasses of morphisms coincide: effective monomorphisms, regular monomorphisms, strong monomorphisms, and extremal monomorphisms. \qed
  \end{corollary}

  Continuing with this line of inquiry, we arrive at the stronger notion of left-invertibility. A morphism in $\barCL$ is left--invertible if and only if it is (the closeness class of) a quasi-isometric embedding whose image is a coarse Lipschitz retract \cite{CH16}.
  We provide some standard examples in Section \ref{sec:control} demonstrating that regular monomorphisms are not necessarily left-invertible.

  \begin{proposition}[RegMono $\neq$ left-invertible]\label{prop:notleft}
   In $\barCL$, the class of regular monomorphisms does not coincide with the class of left-invertible morphisms.
  \end{proposition}

  \begin{remark}
   If $X$ is quasigeodesic and $Y$ is Gromov hyperbolic, then any quasi-isometric embedding $f \from X \to Y$ has quasiconvex image. Thus $f(X)$ is a coarsely Lipschitz retract, hence $f$ admits a left-inverse (up to closeness).
  \end{remark}

  Coregular categories admit a notion of equivalence corelation on objects.
  A coregular category is \emph{coexact in the sense of Barr} if every equivalence corelation arises as a cokernel pair (see Section \ref{sec:coexact}).
  This class of categories fits in between the abelian categories and the coregular categories. The category $\Top$ is an example of a non-coexact coregular category. We show that $\barCL$ also falls into this class.
  \begin{proposition}[Non-coexactness]\label{prop:non-coexact}
   The category $\barCL$ is not coexact in the sense of Barr.
  \end{proposition}

  The results in this paper also hold beyond the coarsely Lipschitz category. We can consider the more general notion of \emph{coarse equivalence} of metric spaces by working in the metric coarse category $\barBorn$ (see Section \ref{sec:control}). We will treat both categories in parallel throughout this paper as most of the arguments will work in either context. Consequently, we obtain analogous results with a few exceptions; we shall remark on these when they arise.

  \begin{theorem}\label{thm:analogy}
   Theorem \ref{thm:Coregular}, Corollaries \ref{cor:Orthogonal} and \ref{cor:Mono}, and Propositions \ref{prop:notleft} and \ref{prop:non-coexact} all hold with $\barBorn$ in place of $\barCL$. Every monomorphism in $\barBorn$ is a regular monomorphism.
  \end{theorem}

  In summary, we obtain the following implications between various subclasses of monomorphisms. The first implication is strict for both $\barCL$ and $\barBorn$, whereas the final implication is strict only for $\barCL$. Note that each forward implication holds in any category.

  \begin{center}
  left-invertible $\implies$ effective monomorphism $\iff$ regular monomorphism $\iff$ \\ strong monomorphism $\iff$ extremal monomorphism $\implies$ monomorphism.
  \end{center}
  
  Our results fit into the broader scheme of understanding coarse geometry from a category theoretic perspective \cite{CH16, DZ17, HW19, Zav19, LV23}. The literature in this direction is more developed for the category of coarse structures in the sense of Roe \cite{Roe03}, which are more general than metric structures. Characterisations or constructions from within this context do not necessarily carry over when restricting to the metric setting. As our primary focus is on the coarsely Lipschitz category, we shall ensure that all arguments, constructions, and proofs take place entirely within the metric context.

  \subsection*{Acknowledgements}
  The author is supported by the National Natural Science Foundation of China (NSFC 12101503); the Gusu Innovation and Entrepreneurship Leading Talents Programme (ZXL2022473); and the XJTLU Research Development Fund (RDF-23-01-121). The author thanks Federico Vigolo for providing helpful comments and references, and Adam--Chistiaan van Roosmalen for interesting conversations on (co)regular and (co)exact categories.

  \section{(Affinely) Controlled maps}\label{sec:control}

  Standard background for coarse geometry can be found in \cite{Roe03, CH16, LV23}.
  We consider only (pseudo)metric spaces, rather than the general coarse structures in the sense of Roe.

  A (possibly multi-valued) map $f \from X \to Y$ between (extended) metric spaces is \emph{controlled} if there exists an increasing function $\Phi \from [0,\infty) \to [0, \infty)$, called an \emph{upper control}, such that
  \[\diam_Y\left(f(x) \cup f(x')\right) \leq \Phi(d_X(x,x'))\]
  for all $x,x' \in X$ satisfying $d_X(x,x') < \infty$. In particular, if $\Phi$ is an upper control for $f$ then $f(x)$ has diameter at most $\Phi(0)$ for all $x \in X$.
  A map $f$ is \emph{coarsely Lipschitz} if it has an affine upper control.
  We call $f$ \emph{uniformly metrically proper} if there is a proper increasing function $\Psi \from [0,\infty) \to [0, \infty)$, called a \emph{lower control}, such that
  \[\Psi(d_X(x,x')) \leq d_Y(y, y')\]
  for all $x,x' \in X$ and $y\in f(x), y'\in f(x')$. A map with both an upper and lower control is called a \emph{coarse embedding}. In particular, if $f$ has affine upper and lower controls then it is called a \emph{quasi-isometric embedding}. We say $f \from X \to Y$ is \emph{coarsely surjective} if there exists some $r \geq 0$ such that the metric $r$--neighbourhood of $f(X)$ equals $Y$. A coarsely surjective coarse (resp. quasi-isometric) embedding is called a \emph{coarse equivalence} (resp. \emph{quasi-isometry}).

  Given $\kappa \geq 0$, two maps $f,g \from X \to Y$ are $\kappa$--\emph{close}, denoted $f \close_\kappa g$, if
  $\diam_Y(f(x) \cup g(x)) \leq \kappa$
  for all $x \in X$. We say that $f$ and $g$ are \emph{close} if they are $\kappa$--close for some $\kappa \geq 0$.
  This yields an equivalence relation on maps. Moreover, closeness is preserved by composition.

  \begin{definition}[$\barBorn$, $\barCL$ categories]
  Let $\Born$ be the category with extended metric spaces as objects and controlled maps as morphisms. Let $\CLip$ be the wide subcategory of $\Born$ where morphisms are coarsely Lipschitz maps. The \emph{metric coarse category} $\barBorn$ (resp. \emph{coarsely Lipschitz category} $\barCL$) is the quotient category of $\Born$ (resp. $\CLip$) with extended metric spaces as objects, and  the closeness classes of controlled (resp. coarsely Lipschitz) maps as morphisms.
  \end{definition}

  \textbf{Notation:} We shall write $\cat$ to denote either the category $\Born$ or $\CLip$, and $\barcat$ for the respective quotient category.
  Given a map $f\in\Mor(\cat)$, write $\bar f \in \Mor(\overline{\cat})$ for its closeness class.

  \begin{lemma}[Initial and terminal objects]
  The empty set is initial in both $\cat$ and $\barcat$. Any non-empty bounded space is terminal in $\barcat$. \qed
  \end{lemma}

  \begin{remark}[Terminology]
   The controlled (resp. coarsely Lipschitz) maps defined here are referred to as \emph{coarsely Lipschitz} (resp. \emph{large-scale Lipschitz}) by Cornulier--de la Harpe \cite{CH16}. Leitner--Vigolo \cite{LV23} refer to coarsely Lipschitz maps as \emph{quasi-Lipschitz}.
   In the setting of metric spaces, the controlled maps (resp. coarse embeddings) defined here are equivalent to the \emph{uniformly bornologous} (resp. \emph{coarse}) maps in the sense of Roe \cite{Roe03}.
  \end{remark}

  \begin{remark}[Multi-- vs. single-- valued, pseudometrics]
  Any multi-valued controlled map is close to a single-valued map. Thus, if we instead insist on morphisms in $\cat$ being single-valued, we would still obtain the same quotient category $\barcat$ defined as above.
  We could also allow for pseudometric spaces. Since any pseudometric space is quasi-isometric to a metric space, for example, by taking its largest Hausdorff quotient, the results in this paper would also hold if we permit $\barcat$ to include pseudometrics.
  \end{remark}

  Let us now present some characterisation of morphisms in $\barcat$ as given by Cornulier--de la Harpe \cite{CH16}. We make slight adjustments to allow for extended metrics. In particular, we do not require the assumption that the domain is non-empty in the characterisation of epimorphisms.
  \begin{theorem}[Morphisms in $\barBorn$ \protect{\cite[Prop 3.A.16]{CH16}}]
  Let $f \from X \to Y$ be a controlled map. Then
  \begin{enumerate}
   \item $f$ has a lower control if and only if $\bar f$ is a monomorphism in $\barBorn$;
   \item $f$ is coarsely surjective if and only if $\bar f$ is an epimorphism in $\barBorn$; and
   \item $f$ is a coarse equivalence if and only if $\bar f$ is an isomorphism in $\barBorn$.
  \end{enumerate}
  Moreover, $\bar f$ is an isomorphism if and only if it is a monomorphism and an epimorphism.
  \end{theorem}
  \proof
  This follows verbatim from \protect{\cite[Prop 3.A.16]{CH16}}, with minor adjustments needed to prove the converse directions of (1) and (2) to allow for extended metric spaces.

  Converse of (1): Assume $f$ does not have a lower control. Then there exists some $\kappa \geq 0$ and sequences $(x_n)_{n\geq 0}$ and $(x'_n)_{n \geq 0}$ in $X$ such that $d_X(x_n, x_n') \geq n$ but $d_Y(fx_n, fx'_n) \leq \kappa$ for all $n$ (noting that $d_X(x_n, x_n')$ could be infinite). Define a metric $d$ on $\N$ by declaring $d(m,n) = \infty$ whenever $m \neq n$. Define 1--Lipschitz maps $g,g' \from \N \to X$ by $g(n) = x_n$ and $g'(n) = x'_n$ for all $n \geq 0$.
  Then $g, g'$ are not close, hence $\bar g \neq \bar g'$ in $\barCL$. However, $fg \approx_\kappa fg'$, hence $\bar f\bar g = \bar f\bar g'$. Thus $\bar f$ is not a monomorphism.

  Converse of (2): Suppose $f$ is not coarsely surjective. Define functions $h,h' \from Y \to [0,\infty]$ by $h(y) = 0$ and $h'(y) = \inf_{x \in X} d_Y(x, y)$ for all $y\in Y$ (noting that if $X = \emptyset$ then $h'(y) = \infty$ for all $y\in Y$). Metrise $[0,\infty]$ by taking the Euclidean metric on $[0,\infty)$ and declaring $d(t,\infty) = \infty$ for $t \in [0,\infty)$. Then $h,h'$ are 1--Lipschitz maps satisfying $hf = h'f$ (this holds vacuously if $X$ is empty). Since $f$ is not coarsely surjective, $h'(Y)$ does not lie in the $\kappa$--neighbourhood of $0$ in $[0,\infty]$ for any $\kappa \geq 0$. Consequently, $h$ and $h'$ are not close and so $\bar f$ is not an epimorphism.
  \endproof

  \begin{theorem}[Morphisms in $\barCL$ \protect{\cite[Prop 3.A.22]{CH16}}]\label{thm:MorCL}
  Let $f \from X \to Y$ be a coarsely Lipschitz map. Then
  \begin{enumerate}
   \item $f$ has a lower control if and only if $\bar f$ is a monomorphism in $\barCL$;
   \item $f$ is coarsely surjective if and only if $\bar f$ is an epimorphism in $\barCL$; and
   \item $f$ is a quasi-isometry if and only if $\bar f$ is an isomorphism in $\barCL$.
  \end{enumerate}
  \end{theorem}
  \proof
  Left or right cancellativity persists when passing to a subcategory, hence the forwards implications of (1) and (2) follow immediately from the previous theorem. The converse for (1) and (2) follow using the same proof as above. Statement (3) is a standard result.
  \endproof

  \begin{remark}[Erratum]
  Theorem \ref{thm:MorCL} appeared as Proposition 3.A.22 in \cite{CH16} without proof. The published version of the book incorrectly claimed that $f$ has an affine lower control if and only if $\bar f$ is a monomorphism in $\barCL$. The corrected statement now appears in an erratum \cite{CH16*}.
  (For the Converse of (1) in the non-extended metric setting, replace $(\N, d)$ in the proof by any real sequence $\{w_n\}_{n\geq 0}$ satisfying $w_n > w_{n-1} + \max \{d_X(x_{n-1}, x_n),  d_X(x'_{n-1}, x'_{n})\}$, regarded as a set equipped with the Euclidean metric, then consider the 1--Lipschitz maps $g(w_n) = x_n$ and $g'(w_n) = x'_n$.)
  \end{remark}

  In contrast to $\barBorn$, the category $\barCL$ is not balanced. That is, a morphism which is monic and epic is not necessarily an isomorphism, as the following example shows.

  \begin{example}[Cubes to squares]
  Consider the map
  \[f \from \{n^3 ~|~ n \in \N \} \to \{n^2 ~|~ n \in \N \}, \qquad n^3 \mapsto n^2,\]
  where each space is equipped with the Euclidean metric.
  This is surjective, 1--Lipschitz, and has a lower control, hence it is a coarse equivalence. Thus $\bar f$ is both a monomorphism and an epimorphism in $\barCL$. However, $f$ is not a quasi-isometry, hence $\bar f$ is not an isomorphism.
  \end{example}

  Another example of a monomorphism in $\barCL$ which is not a quasi-isometric embedding is a horocyle $\gamma \from \R \to \Hy^2$ on the hyperbolic plane, equipped with a unit-speed paramererisation.

  \begin{proposition}[Left-invertibility \protect{\cite[Prop 3.A.20]{CH16}}]
   Let $f \from X \to Y$ be a coarsely Lipschitz (resp. controlled) map. Then $\bar f$ is left-invertible in $\barCL$ (resp. $\barBorn$) if and only if $f$ is a quasi-isometric (resp. coarse) embedding whose image is a coarsely Lipschitz (resp. coarse) retract. \qed
  \end{proposition}
  
  \begin{example}[\protect{\cite[Remark 2.4.10]{LV23}}]
  The inclusion $\{n^2 ~|~ n \in \N\} \hookrightarrow \R$ (with the Euclidean metrics) is an isometric embedding, but $\R$ does not coarsely retract onto its image .
  \end{example}

  Thus, we deduce Proposition \ref{prop:notleft} modulo Theorem \ref{thm:regmono}.
  Another example is the inclusion map of any horocycle into the hyperbolic plane, equipped with the induced metric.

  We conclude this section with a useful observation.
  Given a constant $C \geq 0$, we say that a function $F \from [0,\infty) \to [0,\infty)$ is \emph{$C$--coarsely superadditive} if $F(s) + F(t) \leq F(s+t) + C$ for all $s,t\geq 0$. By the following lemma, we may assume that any upper control satisfies this property.

  \begin{lemma}[Coarse superadditivity]\label{lem:superadditivity}
   Let $F \from [0,\infty) \to [0,\infty)$ be an increasing function. Then $F$ is bounded above by an increasing $C$--coarsely superadditive function for some $C \geq 0$.
  \end{lemma}
  \proof
  Let $C = F(1)$. For each integer $n \geq 1$, the function $F$ is bounded on $[0,n]$, hence there exists some $K_n > 0$ such that $F(t) \leq K_nt + C$ for all $0 \leq t \leq n$. We may assume that the $K_n$ form an increasing sequence. Define $F' \from [0,\infty) \to [0,\infty)$ by setting $F'(0) = C$ and $F'(t) = K_nt + C$ for $t \in (n-1, n]$. Then $F'$ is increasing and bounds $F$ from above.   Given $s,t\geq 0$, let $j = \lceil s \rceil$, $k = \lceil t \rceil$, and $l = \lceil s + t \rceil$. Then
  \[F(s) + F(t) = (K_js + C) + (K_kt + C) \leq K_l(s+t) + 2C = F(s+t) + C, \]
  hence $F'$ is $C$--coarsely superadditive.
  \endproof

  \section{Finite colimits}

   In this section, we show that $\barcat$ admits all finite colimits. This result is standard; see \cite{LV23} for arguments in more general settings. We provide proofs for completeness, and to describe an explicit construction of pushouts which will be convenient for later sections.

   Given a collection of metric spaces $\{(X_i, d_i)\}_{i\in I}$, we equip their disjoint union $\coprod_{i \in I} X_i$ with the metric $d$ defined by $d(x,y) = d_i(x,y)$ whenever $x,y \in X_i$; and $d(x,y) = \infty$ otherwise.

  \begin{lemma}[Finite coproducts exist in $\barcat$]\label{lem:coproduct}
   The coproduct of a finite collection of metric spaces in $\barcat$ is realised by their disjoint union, together with the inclusion maps from each space.
  \end{lemma}

  \proof
  Let $\{(X_i, d_i)\}_{i\in I}$ be a finite collection of metric spaces. Suppose $Z$ is a cocone under the discrete diagram whose objects are the $X_i$'s, together with colegs $\bar \lambda_i \from X_i \to Z$. For each $i \in I$, choose a representative $\lambda_i \in \bar\lambda_i$. Then the $\lambda_i$ admit a common upper control $\Phi = \sum \Phi_i$, where $\Phi_i$ is an upper control for $\lambda_i$. Note that all involved controls can be taken to be affine in the case where $\barcat = \barCL$. Therefore, the map $\lambda := \coprod_i \lambda_i \from \coprod_i X_i \to Z$ is upper controlled by $\Phi$. Therefore, each $\bar\lambda_i$ factors through $\bar\lambda$ via the (closeness class of) the inclusion $\iota_i \from X_i \hookrightarrow \coprod_i X_i$.

  To verify uniqueness, suppose $\bar\mu \from \coprod_i X_i \to Z$ is a morphism in $\barcat$ satisfying $\bar\lambda_i = \bar\mu\bar\iota_i$ for each $i$. Upon choosing a representative $\mu\in \bar\mu$, there exists $\kappa_i \geq 0$ such that $\lambda\iota_i = \lambda_i \approx_{\kappa_i} \mu \iota_i$ for each $i$. For any $x \in \coprod_i X_i$, there exists unique $i$ and $x_i \in X_i$ such that $x = \iota_i x_i$. Setting $\kappa = \max_i \kappa_i$, it follows that $\lambda x = \lambda_i x_i \approx_{\kappa} \mu \iota_i x_i = \mu x$ for all $x \in \coprod_i X_i$. Therefore $\bar\lambda = \bar\mu$ as required.
  \endproof

  Next, we describe a construction commonly used in coarse geometry, where geodesic spaces are glued together via edges.   This is used, for example, to construct coned-off Cayley graphs in Farb's definition of relative hyperbolic groups \cite{Far98}, and in the Bestvina--Bromberg--Fujiwara construction of quasitrees of (geodesic) metric spaces \cite{BBF15}.
  Since we are working with general metric spaces, let us describe the construction without the geodesic assumption.
  \begin{definition}[Coarse gluing]
  Let $(X,d)$ be a metric space and suppose $\sim$ is a symmetric binary relation on $X$.
    Let $\Gamma$ be the metric graph with vertex set $X$, with two types of edges given as follows: for each pair of distinct $x,y \in X$, we connect them by
  \begin{itemize}
   \item an \emph{internal edge} of length $d(x,y)$ if $d(x,y) < \infty$, and
   \item a \emph{glued edge} of length 1 if $x \sim y$.
  \end{itemize}
  Let $\tilde d$ be the induced metric on $X$ under the inclusion $X \hookrightarrow \Gamma$, where $\Gamma$ is equipped with the path metric.
    We call $(X,\tilde d)$ the space obtained by \emph{coarsely gluing} $X$ via $\sim$.
  \end{definition}

  While $(X,\tilde d)$ is not a graph itself, it will be convenient for expository purposes to refer to edges and paths in the overlying graph $\Gamma$. Equivalently, the metric $\tilde d$ is characterised by the following maximality condition: for any metric $d'$ on $X$ satisfying $d' \leq d$ and $d'(x,y) \leq 1$ whenever $x \sim y$, we have $d' \leq \tilde d$.

  In our setting, the binary relations determining the coarse gluings will typically arise from maps between metric spaces. Given a map $f \from X \to Y$, we define the space obtained by \emph{gluing $X$ to $Y$ along $f$} to be the disjoint union $X \sqcup Y$ coarsely glued via the relation given by $x \sim f(x)$ for all $x \in X$. More generally, a disjoint union of metric spaces can be coarsely glued via some family of maps between pairs of spaces. We shall use this to construct pushouts in $\barcat$.

  \begin{proposition}[Pushouts exist in $\barcat$]\label{prop:pushout}
   Let $X_0, X_1, X_2$ be metric spaces and $\bar f_i \from X_0 \to X_i$ be morphisms in $\barcat$ for $i = 1,2$. Choose $f_i \in \bar f_i$ for $i = 1,2$, and let $X_1 \coprod_{X_0} X_2$ be the space obtained by coarsely gluing $X_0$ to $X_1$ and $X_2$ along $f_1$ and $f_2$ respectively.
   Then there is a pushout diagram
   \begin{center}
   \begin{tikzcd}[column sep=normal, row sep=large]
    X_0 \arrow[r, "\bar f_1"] \arrow[d, "\bar f_2"'] & X_1 \arrow[d] \\
    X_2 \arrow[r]  & X_1 \coprod_{X_0} X_2
    \end{tikzcd}
    \end{center}
   in $\barcat$, where the colegs are the (closeness classes of the)
      inclusions $X_i\hookrightarrow X_1 \coprod_{X_0} X_2$ for $i=0,1,2$.
  \end{proposition}

  \proof
  Let $Z$ be a metric space and suppose $\bar\lambda_i \from X_i \to Z$ for $i=0,1,2$ are morphisms in $\barcat$ satisfying $\bar\lambda_i\bar f_i = \bar\lambda_0$ for $i=1,2$. Choose representatives $\lambda_i \from X_i \to Z$ for each $\bar\lambda_i$. Then there is some $\kappa\geq 0$ such that $\lambda_i f_i \approx_\kappa \lambda_0$ for $i=1,2$. Without loss of generality, we may choose a common upper control $\Phi$ for all $\lambda_i$ which is increasing, $C$--coarsely super-additive (and affine in the case of $\barCL$) for some $C \geq 0$, and satisfies $\Phi(1) \geq \kappa$.

  Let $\lambda \from X_1 \coprod_{X_0} X_2 \to Z$ be the map defined on underlying sets by
  \[\lambda := \lambda_0 \sqcup \lambda_1 \sqcup \lambda_2 \from X_0 \sqcup X_1 \sqcup X_2 \to Z.\]
    We wish to find an upper control for $\lambda$.
  Let $d$ denote the metric on $X_1 \coprod_{X_0} X_2$.
  Suppose that $x,x' \in X_1 \coprod_{X_0} X_2$ are points satisfying $d(x,x') < \infty$. For any $\epsilon > 0$, there exists a path $P$ with vertices $x = v_0, v_1, \ldots, v_n = x'$ in $X_1 \coprod_{X_0} X_2$ of length $L \leq d(x,x') + \epsilon$. By deleting vertices if necessary, we may assume that $P$ has no consecutive pair of internal edges. In particular, $P$ has at least $\lfloor \frac{n}{2} \rfloor$ glued edges, hence $n \leq 2d(x,x') + 2\epsilon$.
  By taking $\epsilon > 0$ sufficiently small, we can choose $P$ to have $n \leq \lfloor 2d(x,x') \rfloor$ edges. We shall also assume $\epsilon < 1$.

  Let $v_j, v_{j+1}$ be an edge along $P$. We claim that $d_Z(\lambda v_j, \lambda v_{j+1}) \leq \Phi(d(v_j, v_{j+1}) + 1)$. If this is an internal edge in some $X_i$, then $d_i(v_j, v_{j+1}) < d(v_j, v_{j+1}) + 1$ for otherwise we could reduce the length of $P$ by more than $\epsilon$ by replacing the edge $v_j, v_{j+1}$ by another path in $X_1 \coprod_{X_0} X_2$, contradicting the choice of $P$. Therefore
  \[d_Z(\lambda v_j, \lambda v_{j+1}) \leq \Phi(d_i(v_j, v_{j+1})) \leq \Phi(d(v_j, v_{j+1}) + 1).\]
  If $v_j,v_{j+1}$ is a glued edge, then one endpoint lies in $X_0$ and the other is in its image in $X_i$ under $f_i$ for $i=1$ or $2$. Since $\lambda_i f_i \approx_\kappa \lambda_0$, we also deduce
  \[d_Z(\lambda v_j, \lambda v_{j+1}) \leq \kappa \leq \Phi(1) = \Phi(d(v_j, v_{j+1})) \leq \Phi(d(v_j, v_{j+1}) + 1),\]
  yielding the desired claim.
  Therefore, by respectively applying the triangle inequality; the claim; $C$--coarse super-additivity of $\Phi$; the choice of $P$; and the bound on $n$; we deduce
  \begin{align*}
   d_Z(\lambda x, \lambda x') &\leq \sum_{j=0}^{n-1} d_Z(\lambda v_j, \lambda v_{j+1})\\
   &\leq \sum_{j=0}^{n-1} \Phi(d(v_j, v_{j+1})+1)\\
   &\leq \Phi\left(\sum_{j=0}^{n-1} (d(v_j, v_{j+1}) + 1) \right) + C(n-1)\\
   &\leq \Phi\left(d(x,x') + 1 + n\right) + C(n-1)\\
   &\leq \Phi\left(3 d(x,x') + 1\right) + 2Cd(x,x').
  \end{align*}
    Consequently, $\Phi'(t) := \Phi(3t+1) + 2Ct$ serves as an upper control for $\lambda$. In particular, if $\Phi$ is affine then so is $\Phi'$.

  By construction, each $\bar\lambda_i$ factors through $\bar\lambda$ via the respective (closeness class of the) inclusion $X_i \hookrightarrow X_1 \coprod_{X_0} X_2$.
  Uniqueness follows using a similar argument to Lemma \ref{lem:coproduct} with $X_1 \coprod_{X_0} X_2$ in place of $\coprod_i X_i$.
  \endproof

  Since any category with all finite coproducts and pushouts admits all finite colimits, the following is immediate.

  \begin{proposition}[Finite colimits exist in $\barcat$]\label{prop:colimit}
  The category $\barcat$ admits all finite colimits. \qed
  \end{proposition}

  \begin{remark}[Modified coarse gluing]
  Since any pushout in $\barCL$ (resp. $\barBorn$) is unique up to quasi-isometry (resp. coarse equivalence), there is flexibility in the construction of $X_1 \coprod_{X_0} X_2$. We may, for example, apply any combination of the following modifications at the cost of applying a (canonical closeness class of a) quasi-isometry:
  \begin{itemize}
    \item Allow all glued edges to have lengths uniformly bounded above and below;
    \item For each glued edge, include the whole metric interval as part of the space;
    \item Allow the maps $f_i$ to be multi-valued (the fact they are uniformly controlled means image of each point has uniformly bounded diameter); or
    \item Take the underlying set of $X_1 \coprod_{X_0} X_2$ to be $X_1 \sqcup X_2$, then coarsely glue via the relation $f_1(x) \sim f_2(x)$ for all $x \in X_0$.
  \end{itemize}
  \end{remark}

  \section{Equalisers}

  As observed by Leitner--Vigolo \cite{LV23}, the category $\barcat$ does not admit all equalisers. One example of non-existence will appear in the proof of Prop \ref{prop:coexact}. Nevertheless, we show how they can be computed using a suitable filtration when they do exist.

  Let $\bar f, \bar g \from X \rightrightarrows Y$ be a parallel pair in $\barcat$. Given a constant $\kappa \geq 0$, define the \emph{$\kappa$--equaliser} of $f \in \bar f, g\in \bar g$
  to be
  \[\Eq_\kappa(f,g) := \{x \in X ~|~ \diam_Y(f(x)\cup g(x)) \leq \kappa\}\]
  equipped with the induced metric.
  For any $0 \leq \kappa \leq \kappa'$, the inclusion ${\Eq_\kappa(f,g) \hookrightarrow \Eq_{\kappa'}(f,g)}$ is an isometric embedding. These inclusions assemble to form the \emph{equaliser filtration}
  \[\Eq_*(f,g) \from [0,\infty) \to \barcat \]
  of $f$ and $g$. Note that $f$ and $g$ are $\kappa$--close if and only if $\Eq_\kappa(f,g) = X$.

  \begin{lemma}[Closeness and $\kappa$--equalisers]\label{lem:eq_close}
     Let $f,f',g,g' \from X \to Y$ be maps such that $f \approx_{\kappa'} f'$ and $g \approx_{\kappa'} g'$ for some $\kappa' \geq 0$. Then $\Eq_\kappa(f,g) \subseteq \Eq_{\kappa + 2\kappa'}(f',g')$ for all $\kappa \geq 0$.
  \end{lemma}

    \proof
    Let $x \in \Eq_\kappa(f,g)$. Then $f'x \approx_{\kappa'} fx \approx_{\kappa} gx \approx_{\kappa'} g'x$. Hence  $x \in \Eq_{\kappa + 2\kappa'}(f',g')$.
    \endproof
    Thus, we may choose different representatives of $\bar f, \bar g$ at the cost of modifying the constant $\kappa$. Let us fix a choice of representatives $f,g$, and write $\iota_\kappa \from \Eq_\kappa(f,g) \hookrightarrow X$ for the inclusion.

   \begin{lemma}[Factoring through $\kappa$--equalisers]\label{lem:equaliser_filt}
     Let $\bar h \from W \to X$ be a morphism in $\barcat$ equalising ${\bar f, \bar g \from X \rightrightarrows Y}$. Then $\bar h$ factors uniquely through $\bar \iota_\kappa \from \Eq_\kappa(f,g)\to X$ for all sufficiently large $\kappa \geq 0$.
    \end{lemma}
    \proof
    Let $\bar h \from W \to X$ be a morphism satisfying $\bar f \bar h = \bar g \bar h \from W \to Y$. For any $h \in \bar h$, there exists $\kappa \geq 0$ such that $fh \approx_\kappa gh$.
    Therefore, $\diam_Y(f(h(w)) \cup g(h(w))) \leq \kappa$ for all $w \in W$, hence $h(W) \subseteq \Eq_\kappa(f,g)$. Thus, $h$ factors through $\iota_\kappa \from\Eq_\kappa(f,g) \hookrightarrow X$. To verify uniqueness, suppose that $\bar h = \iota_\kappa \bar k = \iota_\kappa \bar k'$ for some morphisms $\bar k, \bar k'$. Since $\iota_\kappa$ is an isometric embedding, $\bar \iota_\kappa$ is a monomorphism in $\barcat$. Therefore $\bar k = \bar k'$.
    \endproof

    We say that the filtration $\Eq_*(f,g)$ \emph{stabilises} in $\barcat$ if (the closeness class of) the inclusions ${\Eq_\kappa(f,g) \hookrightarrow \Eq_{\kappa'}(f,g)}$ are isomorphisms in $\barcat$ for all sufficiently large $\kappa \leq \kappa'$. Appealing to Lemma \ref{lem:eq_close}, stability of $\Eq_*(f,g)$ depends only on the closeness classes $\bar f, \bar g$.

    \begin{proposition}[Equalisers in $\barcat$]\label{prop:equaliser}
     Let ${\bar f, \bar g \from X \rightrightarrows Y}$ be a parallel pair in $\overline \cat$, and let $f\in \bar f$, $g \in \bar g$ be representatives.
     Then the following are equivalent:
      \begin{enumerate}
       \item There exists $\kappa_0 \geq 0$ such that for all $\kappa \geq \kappa_0$, there exists $r \geq 0$ such that $\Eq_\kappa(f,g)$ lies in the (closed) metric $r$--neighbourhood of $\Eq_{\kappa_0}(f,g)$ in $X$;
       \item $\Eq_*(f,g)$ stabilises in $\barcat$;
       \item $\colimit \Eq_*(f,g)$ exists in $\overline{\cat}$; and
       \item The equaliser of $\bar f, \bar g$ exists in $\barcat$.
      \end{enumerate}
      Moreover, if any of the above hold, then the equaliser for $\bar f, \bar g$ is realised by $\bar\iota_\kappa \from \Eq_\kappa(f,g)\to X$ for sufficiently large $\kappa \geq 0$.
    \end{proposition}

    \proof
    (1) $\iff$ (2). The inclusion $\Eq_{\kappa_0}(f,g) \hookrightarrow \Eq_{\kappa}(f,g)$ is an isometric embedding, hence its closeness class is an isomorphism in $\barcat$ if and only if it is coarsely surjective.

    (2) $\implies$ (3). Trivial.

    (3) $\implies$ (4). Assume $\colimit \Eq_*(f,g)$ exists. Let $\bar\lambda_\kappa \from \Eq_\kappa(f,g) \to \colimit \Eq_*(f,g)$ be the limit colegs. Since the $\bar\iota_\kappa \from \Eq_\kappa(f,g) \to X$ assemble to form the legs of a cocone under $\Eq_*(f,g)$, there is a unique morphism $\bar\iota \from \colimit \Eq_*(f,g)\to X$ such that $\bar\iota\bar\lambda_\kappa = \bar\iota_\kappa$ for all $\kappa \geq 0$.
     \begin{center}
     \begin{tikzcd}[column sep=large, row sep=large]
	{\Eq_\kappa(f,g)} \\
	{\colimit \Eq_*(f,g)} & X && Y 		\arrow["{\bar \iota_\kappa}", from=1-1, to=2-2]
	\arrow["{\bar \mu_\kappa}", shift left, bend left = 15, from=1-1, to=2-4]
	\arrow["{\bar \lambda_\kappa}"', from=1-1, to=2-1]
	\arrow["{\bar f}", shift left=1, from=2-2, to=2-4]
	\arrow["{\bar g}"', shift right=1, from=2-2, to=2-4]
	\arrow["{\bar \iota}"', from=2-1, to=2-2]
	\arrow["{\bar \mu}"', shift right=0.8, bend right = 30, from=2-1, to=2-4]
    \end{tikzcd}
    \end{center}
     Invoking the definition of $\Eq_\kappa(f,g)$, we see that $f\iota_\kappa \approx_\kappa g\iota_\kappa$ for all $\kappa \geq 0$, hence
    \[ \bar f \bar\iota \bar\lambda_\kappa = \bar f \bar\iota_\kappa = \bar g \bar\iota_\kappa = \bar g \bar\iota \bar\lambda_\kappa\]
    for all $\kappa \geq 0$. Therefore, the $\bar\mu_\kappa := \bar f \bar\iota \bar\lambda_\kappa  = \bar g \bar\iota \bar\lambda_\kappa  \from \Eq_\kappa(f,g) \to Y$ assemble to form a cocone under $\Eq_*(f,g)$, and so there exists a unique morphism $\bar \mu \from \colimit \Eq_*(f,g)\to Y$ satisfying ${\bar\mu\bar\lambda_\kappa = \bar\mu_\kappa}$ for all $\kappa \geq 0$.
    This condition is satisfied by both $\bar f \bar\iota$ and $\bar g \bar\iota$, hence $\bar\mu = \bar f \bar\iota = \bar g \bar\iota$.
    In particular, $\bar \iota$ equalises the pair $\bar f, \bar g$.
    Since every cone above $\bar f, \bar g$ factors uniquely through some $\bar\iota_\kappa \from \Eq_\kappa(f,g) \to X$, and every $\bar\iota_\kappa$ factors uniquely through $\bar\iota$, it follows that $\bar\iota \from \colimit \Eq_*(f,g)\to X$ is the desired equaliser.

    (4) $\implies$ (2).
    Assume $\bar h \from E \to X$ is the equaliser of $\bar f, \bar g$. Then, by Lemma \ref{lem:equaliser_filt}, there exists $\kappa_0 \geq 0$ such that $\bar h$ factors uniquely through $\bar \iota_\kappa$ for all $\kappa \geq \kappa_0$.
    Consequently, $\bar \iota_\kappa \from \Eq_\kappa(f,g)\to X$ is the terminal cone above $\bar f, \bar g$ for all $\kappa \geq \kappa_0$. Therefore, the (closeness class of the) inclusion $\Eq_{\kappa_0}(f,g) \hookrightarrow \Eq_\kappa(f,g)$ is an isomorphism in $\barcat$ for all $\kappa \geq \kappa_0$.
    \endproof

  \section{Regular images}

  Since $\barcat$ is not a concrete category, one cannot expect to define the image of a morphism in the conventional sense. Nevertheless, we show that every closeness class admits a regular image, defined in terms of a universal property, which is realised by the (conventional) image of any choice of representative.

  \begin{definition}[Regular image]
  In any category, the \emph{cokernel pair} of a morphism $f \from X \to Y$ is the canonical pair $i_1, i_2 \from Y \rightrightarrows Y \coprod_X Y$ arising from the pushout of the pair $f,f \from X \rightrightarrows Y$. The \emph{regular image} of $f$ is defined as the equaliser
  of its cokernel pair.   \end{definition}

  A morphism which is its own regular image is also called an \emph{effective monomorphism}. Any effective monomorphism is necessarily regular.

  Returning our attention to the category $\barcat$, we may, by Proposition \ref{prop:pushout}, construct the pushout of any morphism $\bar f \from X \to Y$ along itself by choosing any representative $f \in \bar f$, then coarsely gluing $X$ to two disjoint copies of $Y$ along $f$. The cokernel pair $\bar i_1, \bar i_2 \from Y \rightrightarrows Y \coprod_X Y$ is realised by the inclusions $i_1, i_2 \from Y \rightrightarrows Y \coprod_X Y$ to the two respective copies. Write $i_0$ for the inclusion $X \hookrightarrow Y \coprod_X Y$. We show that the regular image of $\bar f$ is realised by the image of $f$ in the usual (concrete) sense.

  \begin{proposition}[Regular images exist in $\barcat$]\label{prop:image}
  The regular image of a morphism $\bar f \from X \to Y$ in $\barcat$ is realised by the inclusion
  \[f(X) \hookrightarrow Y \rightrightarrows Y \coprod_X Y\]
  for any representative $f \in \bar f$, where $f(X)$ is equipped with the induced metric.
  Consequently, every morphism $\bar f \from X\to Y$ in $\barcat$ admits a factorisation $\bar f = \bar m \bar e$, where $\bar m$ is its regular image and $\bar e$ is an epimorphism.
  \end{proposition}

  \proof
  Write $d$ for the metric on $Y \coprod_X Y$.
  Let $x \in X$. Then the path $i_1fx, i_0x, i_2fx$ is formed by two glued edges, hence $d(i_1 fx, i_2 fx) \leq 2$. Since there are no edges between $i_1(Y)$ and $i_2(Y)$, any path from $i_1fx$ and $i_2fx$ must pass through $i_0(X)$, and hence runs across at least two glued edges. Therefore, $d(i_1 fx, i_2 fx) = 2$, and so $f(X) \subseteq \Eq_2(i_1, i_2)$.

  Now suppose $y \in \Eq_\kappa(i_1, i_2)$ for $\kappa \geq 2$. Then $d(i_1 y, i_2 y) \leq \kappa < \infty$ and so there exists a path from $i_1 y$ to $i_2 y$. Any such path must pass through $i_0(X)$, using at least two glued edges, hence
  \[d_Y(y, f(X)) + 2 + d_Y(f(X), y) \leq d(i_1 y, i_2 y). \]
  Therefore,
  \[d_Y(y, f(X)) \quad \leq \quad \frac{d(i_1 y, i_2 y)}{2} - 1 \quad \leq \quad \frac{\kappa}{2} - 1,\]
  hence $\Eq_\kappa(i_1, i_2)$ lies in the closed $(\frac{\kappa}{2} - 1)$--neighbourhood of $f(X) \subseteq \Eq_2(i_1, i_2)$ in $Y$. In particular, by taking $\kappa = 2$, we see that $f(X) = \Eq_2(i_1, i_2)$.
  Moreover, $\Eq_*(i_1, i_2)$ stabilises in $\barcat$ for $\kappa \geq 2$ and so, by Proposition \ref{prop:equaliser}, the inclusion $f(X) \hookrightarrow Y$ represents the desired equaliser.
  \endproof

  We are now ready to prove Theorem \ref{thm:RegMono}.

  \begin{theorem}[Characterisation of regular monomorphisms]\label{thm:regmono}
   Let $\bar f\from X \to Y$ be a morphism in $\barCL$ (resp. $\barBorn$). Then the following are equivalent.
   \begin{enumerate}
    \item $\bar f$ is a regular monomorphism;
    \item any (hence every) $f \in \bar f$ is a quasi-isometric (resp. coarse) embedding; and
    \item the epimorphism $\bar e$ in the factorisation $\bar f = \bar m \bar e$ from Proposition \ref{prop:image} is an isomorphism.
   \end{enumerate}
  \end{theorem}

  \proof
   (1) $\implies$ (2). Suppose $\bar f$ is a regular monomorphism.    By Proposition \ref{prop:equaliser}, we may write $\bar f = \bar \iota \bar u$ where $\bar u$ is an isomorphism in $\barCL$ (resp. $\barBorn$) and $\bar\iota$ is represented by an isometric embedding into $X$. Therefore, up to closeness, any $f\in \bar f$ is an isometric embedding post-composed with a quasi-isometry (resp. coarse equivalence), and is hence a quasi-isometric (resp. coarse) embedding.

   (2) $\iff$ (3).  Let $f\in\bar f$.
    By Proposition \ref{prop:image}, we may factorise $\bar f = \bar m \bar e$ via the closeness classes represented by the inclusion $m \from f(X) \hookrightarrow Y$ and corestriction $e \from X \to f(X)$ of $f$. Since $e$ is surjective, the map $f \from X \to Y$ is a quasi-isometric (resp. coarse) embedding if and only if $e$ is a quasi-isometry (resp. coarse equivalence). This holds precisely when $\bar e$ is an isomorphism in $\barCL$ (resp. $\barBorn$).

    (3) $\implies$ (1).    If $\bar e$ is an isomorphism, then $\bar f$ equalises the same parallel pairs as $\bar m$. Since $\bar m$ equalises the cokernel pair of $\bar f$, it follows that $\bar f$ is a regular monomorphism.
     \endproof

  \begin{corollary}[Monomorpisms in $\barBorn$]\label{cor:coarse_mono}
   Every monomorphism in $\barBorn$ is regular. \qed
  \end{corollary}

  This establishes the second sentence of Theorem \ref{thm:analogy}.

  \section{Coregularity and orthogonal factorisation} \label{sec:coregular}

  In this section, we establish coregularity of the category $\barcat$, hence proving Theorem \ref{thm:Coregular}. We then apply standard results to deduce uniqueness of the factorisation from Proposition \ref{prop:image}, and give several equivalent characterisations of quasi-isometric (resp. coarse) embeddings. Background on (co)regular categories can be found in \cite{BGO71, Bor94, Gra21}.

  \begin{theorem}[$\barcat$ is coregular]
  The category $\barcat$ is \emph{coregular}. That is, $\barcat$ admits all finite colimits and regular images, and regular monomorphisms are stable under pushouts.
   \end{theorem}
  \proof
  The existence of finite colimits and regular images directly follow from Propositions \ref{prop:colimit} and \ref{prop:image} respectively.

  To verify the final condition, suppose $\bar m \from X \to Y$ is a regular monomorphism and $\bar f \from X \to Z$ is any morphism. Appealing to Proposition \ref{prop:image} and Theorem \ref{thm:regmono}, we may assume, without loss of generality, that $\bar m$ is represented by an isometric embedding $m \from X \hookrightarrow Y$.
  Consider the following commutative diagram, where the upper and outer squares are pushouts. We wish to prove that the pushout $\bar f_* \bar m$ of $\bar m$ along $\bar f$ is a regular monomorphism.
    \begin{center}
        \begin{tikzcd}[column sep=large]
    X \arrow[rrrr, "\bar f"] \arrow[dd, "\bar m"', tail] \arrow[rd, "\bar 1_X"', "\cong"] & & & & Z \arrow[ld, "\bar f_* \bar 1_X", "\cong"'] \arrow[dd, "\bar f_* \bar m"] \\                                                                                                                         & X \arrow[ld, "\bar m"', tail] \arrow[rr, "(\bar 1_X)_*\bar f"'] & & X \coprod_X Z \arrow[rd, "\bar g", dashed] &                                                                                     \\
    Y \arrow[rrrr, "\bar m_* \bar f"]                                                             & & & & Y\coprod_X Z
    \end{tikzcd}
    \end{center}
    By the universal property, there is a unique morphism
    $\bar g\from X \coprod_X Z \to Y \coprod_X Z$ making the diagram commute. Since $\bar 1_X$ is an isomorphism, its pushout $\bar f_* \bar 1_X$ along $\bar f$ is also an isomorphism. Therefore, $\bar f_* \bar m$ is a regular monomorphism if and only if the same is true for $\bar g$.

    By Proposition \ref{prop:pushout}, we may realise $X \coprod_X Z$ (resp. $Y \coprod_X Z$) by coarsely gluing $X$ to $X$ (resp. $Y$) and $Z$ along $1_X$ (resp. $m$) and $f$.
    Consequently, we may represent $\bar g$ using the map given (on underlying sets) by
    \[g := 1_X \sqcup m \sqcup 1_Z \from X \sqcup X \sqcup Z \to X \sqcup Y \sqcup Z.\]
    We claim that $g$ is an isometric embedding. By Theorem \ref{thm:regmono}, it would follow that $\bar g$ is a regular monomorphism.
    Since $1_X$, $m$, and $1_Z$ are isometric embeddings, $g$ maps each internal edge in $X \coprod_X Z $ to an internal edge of the same length in $Y \coprod_X Z $. The map $g$ also preserves the length of glued edges since, by construction, they all have unit length. Therefore, $g$ preserves the lengths of paths, and is hence 1--Lipschitz.

    Now let $w,w' \in g(X \coprod_X Z) \subseteq Y \coprod_X Z$, and suppose there exists a path ${w = v_0, \ldots, v_n = w'}$ in $Y \coprod_X Z$. If some vertex $v_j$ does not lie in $g(X \coprod_X Z)$, then by deleting this vertex from the sequence we obtain a path which is no longer than the original path. Therefore, the distance between $w,w'$ in $Y \coprod_X Z$ is the infimal length among all paths between them that lie in $g(X \coprod_X Z)$. Each path in $g(X \coprod_X Z)$ is the image of some path in $X \coprod_X Z$ under $g$ (in fact, such a path is unique since $g$ is injective). Since $g$ preserves the lengths of paths, it follows that $g$ cannot decrease distance. Therefore, $g$ is an isometric embedding.
    \endproof

    One feature of coregular categories is that they admit (Epi, RegMono)--orthogonal factorisation systems; see \cite{Bou77, FK72, AHS06} for further background.
    In any category, a morphism $e$ is \emph{left-orthogonal} to a morphism $m$ (and $m$ is \emph{right-orthogonal} to $e$) if given any commutative square of the form
   
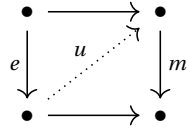
\begin{center}
    \begin{tikzcd}[column sep=large, row sep=large]
    \bullet \arrow[r] \arrow[d, "e"'] & \bullet \arrow[d, "m"] \\
    \bullet \arrow[r] \arrow[ru, "u", dotted] & \bullet
    \end{tikzcd}
    \end{center}
    there exists a unique morphism $u$, called the \emph{diagonal filler}, making the diagram commute.
    If $\mathcal{N}$ is a class of morphisms in a category, write $^\perp\mathcal{N}$ (resp. $\mathcal{N}^\perp$) for the class of morphisms that are left-- (resp. right--) orthogonal to all morphisms in $\mathcal{N}$.

  \begin{definition}[Orthogonal factorisation system]\label{def:OFS}
   Let $\mathcal{E}$, $\mathcal{M}$ be subclasses of morphisms in any category. We say that the category admits an $(\mathcal{E}, \mathcal{M})$--\emph{orthogonal factorisation system} if every morphism $f$ factors as $f = me$ for some $e \in \mathcal{E}$ and $m \in \mathcal{M}$, and either of the following equivalent conditions hold:
   \begin{enumerate}
    \item $\mathcal{E}^\perp = \mathcal{M}$ and $^\perp\mathcal{M} = \mathcal{E}$; or
    \item the factorisation is unique up to unique isomorphism, and both $\mathcal{E}$, $\mathcal{M}$ are closed under composition and contain all isomorphisms.
      \end{enumerate}
   In particular, each of the subclasses $\mathcal{E}$ and $\mathcal{M}$ uniquely determines the other.
  \end{definition}

  \begin{proposition}[(Epi, RegMono)--factorisation]\label{prop:factor}
   Let $\mathcal{E}$, $\mathcal{M}$ respectively be the subclasses of epimorphisms and regular monomorphisms in a coregular category.
   Then the category admits an $(\mathcal{E}, \mathcal{M})$--orthogonal factorisation system. \qed
  \end{proposition}

  \begin{corollary}[Unique factorisation]
   The category $\barcat$ admits an (Epi, RegMono)--orthogonal factorisation system. Moreover, suppose $\bar f = \bar m \bar e$ is the factorisation of a morphism in $\barcat$ as given in Proposition \ref{prop:image}. Let $\bar f = \bar m'\bar e'$ be any factorisation via an epimorphism $\bar e'$ and regular monomorphism $\bar m'$. Then there is a unique morphism $\bar u$ such that $\bar e' = \bar u \bar e$ and $\bar m = \bar m' \bar u$. \qed
  \end{corollary}

  This yields Corollary \ref{cor:Orthogonal}.
  Finally, we invoke the equivalence of several classes of monomorphisms in coregular categories to obtain further characterisations of quasi-isometric (resp. coarse) embeddings.

  \begin{proposition}[Equivalent classes of monomorphisms]\label{prop:equiv_mono}
   Let $m$ be a morphism in a coregular category. Then the following are equivalent:
   \begin{enumerate}
    \item $m$ is an \emph{effective monomorphism}: $m$ is the equaliser of its cokernel pair;
     \item $m$ is \emph{regular monomorphism}: $m$ is the equaliser of some parallel pair;
    \item $m$ is a \emph{strong monomorphism}: $m$ is right-orthogonal to every epimorphism;
    \item $m$ is an \emph{extremal monomorphism}: if $m = ge$ for some epimorphism $e$, then $e$ is an isomorphism. \qed
   \end{enumerate}
    \end{proposition}

  Corollary \ref{cor:Mono} then follows from Theorem \ref{thm:regmono} and Proposition \ref{prop:equiv_mono}. By Corollary \ref{cor:coarse_mono}, the class of monomorphisms in $\barBorn$ coincides with each class given in Proposition \ref{prop:equiv_mono}.

  An explicit construction of diagonal fillers in $\barcat$ is given in Appendix \ref{app:diagonal}.

  \section{Non-coexactness}\label{sec:coexact}

  In this section, we show that $\barcat$ is not coexact in the sense of Barr. See \cite{BGO71, Bor94} for background on (co)exact categories.

  \begin{definition}[Corelation] In a category with finite coproducts, a \emph{corelation} on an object $X$ is a pair of morphisms $f,g \from X \to R$ such that the canonical morphism $f \coprod g \from X \coprod X \to R$ is an epimorphism.
   \end{definition}

  \begin{definition}[Equivalence corelation]
   A corelation $f^+, f^- \from X \to R$ on an object $X$ in a category with finite colimits is an \emph{equivalence corelation} if there exists:
   \begin{enumerate}
    \item a reflexivity morphism $\rho \from R \to X$ satisfying $\rho f^\pm = 1_X$;
    \item a symmetry morphism $\sigma \from R \to R$ satisfying $\sigma f^\pm = f^\mp$; and
    \item a transitivity morphism $\tau \from R \to R \coprod_X R$     satisfying $\tau f^\pm = g^\pm f^\pm$, where $g^\pm$ are as in the pushout diagram:
    \begin{center}
        \begin{tikzcd}[column sep=normal, row sep=large]
    X \arrow[r, "f^+"] \arrow[d, "f^-"'] & R \arrow[d, "g^-"] \\
    R \arrow[r, "g^+"']  & R \coprod_X R
    \end{tikzcd}
    \end{center}
   \end{enumerate}
  The morphisms $\rho$, $\sigma$, and $\tau$ are unique if they exist.
  \end{definition}

  If $f \from X \to Y$ is any morphism then its cokernel pair $Y \rightrightarrows Y \coprod_X Y$, if it exists, yields an equivalence corelation on $Y$. In a coexact category, these account for all the equivalence corelations.

  \begin{definition}[Coexact category]
   An equivalence corelation is \emph{effective} if its equaliser exists, and it is the cokernel pair of its equaliser.
   A coregular category is \emph{coexact in the sense of Barr} if every equivalence corelation is effective.
  \end{definition}

    We conclude with a proof of Proposition \ref{prop:non-coexact}.

  \begin{proposition}[$\barcat$ is not coexact]\label{prop:coexact}
   The category $\barcat$ is not coexact in the sense of Barr.
  \end{proposition}

  \proof
    Our goal is to construct a non-effective equivalence corelation.
    For each integer $n \geq 0$, let $X_n = [1,\infty)$ equipped with Euclidean metric, and $Y_n = \R^2$ equipped with the $L^1$--metric. Let $X$ be the disjoint union of all $X_n$, with $1\in X_m$ coarsely glued to $1\in X_n$ for all $m\neq n$. Similarly, let $Y$ be the disjoint union of all the $Y_n$, with $(0,0)\in Y_m$ coarsely glued to $(0,0)\in Y_n$ for all $m\neq n$.
    For notational convenience, we shall identify $X$ with $\N \times [1,\infty)$ and $Y$ with $\N \times \R^2$.     With this identification, the metrics are given explicitly by
    \begin{align*}
     d_X((m,s), (n,t)) &= \delta_{mn}|t-s| + (1 - \delta_{mn})(t + s - 1) \\
     d_Y((m,x,y), (n,x',y')) &= \delta_{mn}(|x-x'| + |y-y'|) + (1 - \delta_{mn})(|x| +|x'| + |y| + |y'| + 1),
    \end{align*}
    where $\delta_{mn}$ is the Kronecker delta.
    For each $n \geq 0$, let $f_n^\pm \from X_n \to Y_n$ be the (unit-speed geodesic) path given by
    \[f^\pm_n(t) = \begin{cases}
                    (0, \pm t), & 1\leq t \leq n \\ (t-n, \pm n), & t \geq n.
                   \end{cases}
    \]
    Define $f^\pm \from X \to Y$ by $f(n,t) = (n, f^\pm_n(t))$; this restricts to $f^\pm_n$ on each $X_n$.     Replace the codomain of $f^\pm$ by the joint image $R = f^+(X) \sqcup f^-(X) \subset Y$ equipped with the induced metric. Thus, $f^\pm\from X \to R$ are jointly surjective, hence $\bar f^+, \bar f^-\from X \to R$ is a corelation on $X$.
    We claim that $\bar f^+, \bar f^-\from X \to R$ is an equivalence corelation.
    Define maps $\rho \from R \to X$ and $\sigma \from R \to R$ by
    \[\rho(n,x,y) = (n, |x| + |y|) \qquad \textrm{and} \qquad \sigma(n,x,y) = (n,x,-y)\]
    respectively. Note that $\rho$ is 1--Lipschitz and $\sigma$ is an isometry. These maps satisfy $\rho f^\pm = 1_X$ and $\sigma f^\pm = f^\mp$, hence $\bar \rho$ and $\bar \sigma$ serve as the required reflexivity and symmetry morphisms.
    By Proposition \ref{prop:pushout}, we may construct the pushout $R \coprod_X R$ by coarsely gluing $X$ to two disjoint copies $R^+$, $R^-$ of
    $R$ via $f^+$ and $f^-$ respectively,   with the colegs realised by inclusion maps. Denote these inclusions by $g^\pm \from R \hookrightarrow R \coprod_X R$ (mapping to the copy $R^\pm$) and $\iota \from X \hookrightarrow R \coprod_X R$. Thus, the coarse gluing relation is given by $g^+f^+x \sim \iota x \sim g^-f^- x$ for all $x \in X$.     Consequently, the (1--Lipschitz) map $\tau := g^+f^+\rho \from R \to R \coprod_X R$ satisfies
    \[\tau f^\pm = g^+f^+\rho f^\pm = g^+f^+ \approx_2 g^- f^-. \]
    Therefore, $\bar \tau$ serves as the required transitivity morphism.

    Finally, we compute the equaliser filtration of the pair $f^+, f^-$.
    Observe that
    \[d_Y(f^+(n,t), f^-(n,t)) = \min \{2n,2t\}.\]
    Thus, for any $\kappa \geq 0$, we have $(n,t) \in \Eq_{2\kappa}(f^+, f^-) \subseteq X$ if and only if $\min \{n,t\} \leq \kappa$. Therefore,
    \[\Eq_{2\kappa}(f^+, f^-) = \left(\N \cap [1,\kappa]\right) \times [1,\infty) \; \bigsqcup \;  (\N \cap (\kappa, \infty)) \times [1,\kappa]. \]
    Consequently, $\Eq_{2n}(f^+, f^-) \hookrightarrow \Eq_{2n+2}(f^+, f^-)$ is not coarsely surjective for all $n \geq 0$, and so $\Eq_*(f^+, f^-)$ does not stabilise in $\barcat$. By Proposition \ref{prop:equaliser}, the equaliser of $\bar f^+, \bar f^- \from X \to R$ does not exist in $\barcat$, and so $\bar f^+, \bar f^-\from X \to R$ is a non-effective equivalence corelation on $X$.   \endproof

  \begin{remark}
   Strictly speaking, the coarse gluing step in the construction of $X$ and $Y$ is not required for the proof. This is done to demonstrate the existence of a non-effective equivalence corelation on a non-extended metric space.
  \end{remark}

  \appendix
  \section{Diagonal fillers}\label{app:diagonal}

  In this appendix, we give an explicit construction of the diagonal fillers as guaranteed by the right-orthogonality of the regular monomorphisms to the epimorphisms in $\barcat$.
  Let us recall some elementary observations concerning binary relations between metric spaces.

  Regard a binary relation $f \from X \to Y$ between sets $X,Y$ as a subset $f \subseteq X\times Y$. The transpose $f^T \from Y \to X$ of $f$ is defined by declaring $(y,x) \in f^T \iff (x,y) \in f$. The composite of $f \from X \to Y$ and $g \from Y \to Z$ is defined by $(x,z) \in gf \iff \exists y\in Y$ such that $(x,y) \in f$ and $(y,z) \in g$. Note that $(fg)^T = g^Tf^T$. A relation $f \from X \to Y$ is total (i.e. $\forall x\in X,~ \exists y\in Y$ such that $(x,y) \in f$) if and only if $f^T$ is surjective; this holds precisely when $f^Tf$ contains the diagonal $\Delta_X \subset X \times X$. If $f \subseteq f' \from X \to Y$ and $g \subseteq g' \from Y \to Z$ then $gf \subseteq g'f'$.

  Given a metric space $(X,d_X)$, the (closed) $t$--neighbourhood relation $N_t\subseteq X \times X$ is defined by $(x,x')\in N_t \iff d_X(x,x') \leq t$. Observe that $N_sN_t \subseteq N_{s+t}$ and $(N_t)^T = N_t$. If a binary relation $f \from (X,d_X) \to (Y,d_Y)$ between metric spaces has upper control $\Phi$ then $fN_t \subseteq N_{\Phi (t)} f$ and $ff^T \subseteq N_{\Phi(0)}$. A relation $f$ has lower control $\Psi$ if and only if its transpose $f^T$ has upper control $\Psi^T(t) := \sup \{s \geq 0 ~|~ \Psi(s) \leq t\}$. Relations $f,f' \from X \to Y$ satisfy $f \approx_\kappa f'$ if and only if $f\subseteq N_\kappa f'$ and $f' \subseteq N_\kappa f$. If $f$ has upper control $\Phi$ and $f' \subseteq N_\kappa f$ then $f'$ has upper control $\Phi'(t) = \Phi(t) + \kappa$; in addition, if $f'$ is total then $f \subseteq N_{\kappa + \Phi(0)} f'$.
  \begin{lemma}[Diagonal fillers in $\barcat$]
  Let $f,g,e,m$ be morphisms in $\cat$, as in the following diagram, which admit a common upper control $\Phi$ and satisfy $mf \approx_\kappa ge$ for some $\kappa \geq 0$. Assume that $N_r e(X) = Z$ for some $r \geq 0$ and that $m$ has lower control $\Psi$. We further assume that $\Psi$ is affine in the case $\cat = \CLip$.
     \begin{center}
        \begin{tikzcd}[column sep=large, row sep=large]
    X \arrow[r, "f"] \arrow[d, "e"'] & Y \arrow[d, "m"] \\
    Z \arrow[r, "g"'] \arrow[ru, "u", dotted] & W
    \end{tikzcd}
    \end{center}
   Then the binary relation $u := fe^TN_r \from Z \to Y$ is total, has     upper control
    \[\Phi'(t) := \Psi^T\Phi(t) + \Psi^T(\kappa + \Phi(\Phi(0)+ r)),\]
    and satisfies $g \approx_{\kappa + 2\Phi(0) + \Phi(\Phi(0)+ r)} mu$ and $f \approx_{\Phi'\Phi(0)} ue$.
   In particular, $\bar u$ is the unique diagonal filler for the corresponding commuting square in $\barcat$.
      \end{lemma}

  \proof
   The given morphisms $f,g,e,m$ are total.
   The relation $u$ is total since $(e^TN_r)^T = N_re$ is surjective. Moreover,
      \begin{align*}
   u &= fe^TN_r \subseteq m^Tmfe^TN_{r} \subseteq m^TN_\kappa gee^TN_{r} \subseteq m^TN_\kappa gN_{\Phi(0)}N_{r} \\
   & \subseteq m^TN_{\kappa + \Phi(\Phi(0)+ r)}g \subseteq N_{\Psi^T(\kappa + \Phi(\Phi(0)+ r))}m^Tg.
   \end{align*}
      Therefore, $u$ has upper control $\Phi'(t) := \Psi^T\Phi(t) + \Psi^T(\kappa + \Phi(\Phi(0)+ r))$. If $\Phi, \Psi$ are both affine, then so is $\Phi'$. Thus $\bar u$ is a morphism in $\barcat$. Since $mu$ is total and satisfies
   \[mu \subseteq mm^TN_{\kappa + \Phi(\Phi(0)+ r)}g \subseteq N_{\Phi(0)}N_{\kappa + \Phi(\Phi(0)+ r)}g, \]
   we deduce $g \approx_{\kappa + 2\Phi(0) + \Phi(\Phi(0)+ r)} mu$, using the fact that $g$ has upper control $\Phi$.
   Furthermore,
   \begin{align*}
    f \subseteq fe^Te \subseteq fe^TN_re = ue.
   \end{align*}
    Since $f$ is total and $ue$ has upper control $\Phi'\Phi$, it follows that $f \approx_{\Phi'\Phi(0)} ue$.
            Therefore, $\bar u$ serves as a desired diagonal filler in $\barcat$.
    Uniqueness follows from coregularity of $\barcat$ and the fact that $\bar e$ is an epimorphism and $\bar m$ a regular monomorphism.
  \endproof

\providecommand{\bysame}{\leavevmode\hbox to3em{\hrulefill}\thinspace}
\providecommand{\MR}{\relax\ifhmode\unskip\space\fi MR }
\providecommand{\MRhref}[2]{%
  \href{http://www.ams.org/mathscinet-getitem?mr=#1}{#2}
}
\providecommand{\href}[2]{#2}

  \end{document}